\documentclass[12pt,a4paper]{article}

\usepackage{amsmath}
\usepackage{amssymb}
\usepackage{latexsym}
\usepackage{ifthen}
\usepackage{epsfig}

\newcommand{\defi}[1]{\begin{description}\item[Definition:]{#1}\end{description}}
\newcommand{\lkb}{\cdots}
\newcommand{\G}{\Gamma}
\newcommand{\Gh}{\hat{\G}}
\renewcommand{\l}{\lambda}
\newcommand{\q}{\bf{q}}
\newcommand{\ot}{\otimes}

\newcommand{\qbin}[2]{\binom{#1}{#2}_{\hspace{-1ex}q\hspace{0ex}}}
\newcommand{\qibin}[2]{\binom{#1}{#2}_{\hspace{-1ex}q_i\hspace{0ex}}}
\renewcommand{\a}{a_{ij}}
\newcommand{\abild}[1]{\begin{center}\epsfbox{#1}\end{center}}
\newcommand{\pth}[1]{#1^{\text{th}}}
\newcommand{\Z}{\mathbb Z}
\newcommand{\N}{\mathbb N}

\newcommand{\uqg}{U_q(g)}
\newcommand{\D}{\mathcal D}
\newcommand{\ur}{\mathfrak u}

\newcommand{\beer}[3][leer]{\ifthenelse{\equal{#1}{leer}}{\begin{#2}}{\begin{#2}[#1]} \begin{description} \item[{}] #3 \ifthenelse{\equal{#2}{exp}}{\mbox{\hspace{2em}} \hspace{\fill} $\bigtriangleup$}{} \end{description}\end{#2}}
\newcommand{\prf}[2][leer]{\begin{description} \ifthenelse{\equal{#1}{leer}}{\item[Proof:]}{\item[Proof #1:]} #2 \hspace{5em}  \hspace*{\fill}\bf{qed.} \end{description} } 

\DeclareMathOperator{\ad}{ad}
\DeclareMathOperator{\del}{\Delta}
\DeclareMathOperator{\eps}{\varepsilon}
\DeclareMathOperator{\an}{S}

\begin{document}

\author{Daniel Didt\footnote{This work will be part of the author's Ph.D. thesis written under the supervision of Professor H.-J. Schneider. The author is a member of the Graduiertenkolleg ``Mathematik im Bereich ihrer Wechselwirkung mit der Physik'' at Munich University.}\\\small Mathematisches Institut, LMU M\"unchen, Theresienstr. 39\\\small 80333 M\"unchen, Germany\\\small didt@mathematik.uni-muenchen.de}
\title{Linkable Dynkin Diagrams}
\date{}
\maketitle

\begin{abstract}
In this article we develop some aspects of the construction of new Hopf algebras found recently by Andruskiewitsch and Schneider \cite{AS}. There the authors classified (under some slight restrictions) all pointed finite dimensional Hopf algebras with coradical $(\Z/p)^s.$
We contribute to this work by giving a closer description of the possible ``exotic'' linkings.
\end{abstract}

\section{Introduction}
In a series of papers concerned with pointed Hopf algebras, Andruskiewitsch and Schneider developed the ``lifting method''. The application of this idea to finite dimensional Hopf algebras constructed from Dynkin diagrams shows that one obtains a whole class of Hopf algebras by considering ``linkings'' in the Dynkin diagram \cite{AS}.
In this spirit we first introduce a new class of infinite dimensional Hopf algebras that are variations of quantized enveloping algebras of Kac-Moody algebras. 
For this we assign to every Cartan matrix, not necessarily of finite type, a finite number of linkable Dynkin diagrams. For some of those, a so called linkable braiding matrix can be constructed. Together with an abelian group that realizes such a braiding matrix, these Hopf algebras can then be defined.\\
The main part of this article will then be concerned with a detailed investigation of when a linkable braiding matrix of a given Cartan type does exist. This will lead to a characterization of the corresponding linkable Dynkin diagrams. We show how these ideas are related to the usual quantized enveloping algebras and to the finite dimensional Hopf algebras constructed in \cite{AS}, which are themselves variations of the Frobenius-Lusztig kernels \cite{FL}.\\
Another application is the existence of such ``exotic'' linkings as an even number of copies of diagrams $A_n,\,n>2,$ linked into a circle.\\ 
With the structure theorem for linkable Dynkin diagrams, the first main step towards the explicit description of all liftings of Nichols algebras for a given type of diagram has been taken. The question of which groups can realize a given linkable braiding matrix however, must still be addressed. At the end we will discuss an aspect of this for the groups $(\Z/p)^2.$\\

To simplify some technical arguments, we require the base field $k$ to contain a $\pth p$ root of unity for some prime $p>3.$\\
For $q\in k$ we denote the q-binomial as usual by$$\qbin ni=\frac{(n)!_q}{(n-i)!_q(i)!_q}\quad,$$ where $(n)!_q=(n)_q\ldots(2)_q(1)_q,$ and $(n)_q=1+q+\ldots+q^{n-1}$ for $n\in\N,\; 0\leq i\leq n.$ For reasons of completeness we list a few important identities, which can be shown by direct calculation.
\begin{align}
\label{qident1}q^i\qbin ni+\qbin n{i-1}=\qbin ni+q^{n+1-i}\qbin n{i-1}&=\qbin{n+1}i,\quad 1\leq i\leq n\\
\label{qident2}\sum_{i=0}^n(-1)^iq^{\binom i2}\qbin ni&=0,\quad n\geq 1\\
\label{qident3}\sum_{i=0}^n(-1)^iq^{(i^2+i)/2-ni}\qbin ni&=0,\quad n\geq 1
\end{align}

\section{A class of Hopf algebras}
\subsection{Construction}
Let $(\a)$ be a generalized $(s\times s)$-Cartan matrix (cf. \cite{KAC}). The corresponding Dynkin diagram with a number of additional edges, drawn as dotted edges that do not share vertices, will be denoted $D$ and called a \emph{linkable Dynkin diagram}. Two vertices $i$ and $j\neq i$  connected by such dotted edges are called \emph{linkable}. This is written $i\lkb j.$
\defi{A linkable braiding matrix of $D-$Cartan type is an $(s\times s)$ matrix $(b_{ij})$ with the following properties
\begin{align}
b_{ii}&\neq 1\\
\label{cartan} b_{ij}b_{ji}&=b_{ii}^{\a}\\
\label{link}  b_{ki}^{1-\a}b_{kj}&=1 \quad k=1,\ldots,s\,\text{ if }i \text{ is linkable to } j .
\end{align}
Such a matrix is called \emph{realizable over the abelian group $\G$} if there are elements $g_1,\ldots,g_s\in\G$ and characters $\chi_1,\ldots\chi_s\in\Gh$ such that 
\begin{align} 
b_{ij}=\chi_j(g_i)\qquad\text{ for all }i,j\quad\text{ and}\\
\chi_i^{1-\a}\chi_j=1\qquad\text{ whenever }i\lkb j.
\end{align}}
For the construction of our Hopf algebras we need the following ingredients: a linkable Dynkin diagram $D$, a corresponding linkable braiding matrix $(b_{ij})$ of $D$-Cartan type, an abelian group $\G$ over which $(b_{ij})$ can be realized along with the elements $g_i$ and characters $\chi_i$. Furthermore, we choose elements $h_1,\ldots h_l\in\G,$ such that $\G=<h_1>\oplus<h_2>\oplus\ldots\oplus<h_l>.$ Finally we need a collection of numbers $(\l_{ij}\in\{0,1\})_{1\le i<j\le s}$ where $\l_{ij}=0$ if $i$ is not linkable to $j.$ Vertices $i$ and $j$ with $\l_{ij}=1$ are then called \emph{linked}. The collection of objects described above will be called a \emph{linking Datum of finite Cartan type for $\G$}.
\defi{For a linking Datum $\D$ of finite Cartan type for $\G$ we define an algebra $\ur(\D)$ generated by elements $h_1,\ldots,h_l$, $a_1,\ldots,a_s$ with the group relations from $\G$ among the $h_i$'s and the defining relations
\begin{align}
\label{mixed}h_ia_j&=\chi_j(h_i)a_jh_i\qquad 1\leq i\leq l,\, 1\leq j\leq s,\\
\label{serre}(\ad a_i)^{1-\a}(a_j)&=\l_{ij}(1-g_i^{1-\a}g_j)\qquad 1\leq i<j\leq s.
\end{align}
In this case the adjoint representation $\ad$ has the following explicit form:
\begin{equation}
(\ad a_i)^{1-\a}(a_j):=\sum_{k=0}^{1-\a}(-1)^k\qibin{1-\a}kq_i^{\binom k2}b_{ij}^ka_i^{1-\a-k}a_ja_i^k,
\end{equation}
where $q_i:=b_{ii}=\chi_i(g_i).$}
\beer{propo}{
The so defined algebra $\ur(\D)$ has the structure of a Hopf algebra determined by the comultiplication on the generators
\begin{equation}
\del(h_i):=h_i\ot h_i\qquad\del(a_i):=a_i\ot 1+g_i\ot a_i.
\end{equation}}

For the proof of this statement one has to check that the imposed relations define a Hopf ideal. This is a basic exercise. Only the ``quantum Serre'' relations (\ref{serre}) need some extra attention. But both sides of (\ref{serre}) are $(g_i^{1-\a}g_j,1)$-primitive. For the left hand side one can use for instance \cite[Lemma A.1.]{AS-cartan}.

\subsection{Connection with ${\bf\uqg}$}
We want to show how these new Hopf algebras are connected with previously known ones. For instance we can reproduce the usual quantized Kac-Moody Hopf algebras $\uqg.$ \\
We start with the direct sum of two copies of the given symmetrisable Cartan matrix. In the associated Dynkin diagram we connect corresponding vertices by dotted lines. We number the vertices of one copy of the original diagram from 1 to $N$ and the remaining ones from $N+1$ to $2N$ in the same order. The group $\G$ is simply $\Z^N.$ We take the canonical basis of $\G$ to be the $g_i,\; 1\leq i\leq N$, set $g_{N+i}=g_i$ and define characters $\chi_j(g_i):=q^{d_i\a}, \chi_{N+i}=\chi_i^{-1},$ where $d_i\a=d_ja_{ji}.$ As a linkable braiding matrix of the given Cartan type we can now take $b_{ij}=\chi_j(g_i).$ We set $\l_{i(N+i)}=1,\; 1\leq i\leq N,$ and all others 0. The Hopf algebra obtained from this complete linking datum by the above construction is the quantized Kac-Moody algebra. To see this, one sets $K_i:=g_i, K^{-1}_i:=g_i^{-1},E_i:=a_i,F_i:=(q^{-d_i}-q^{d_i})^{-1}a_{N+i}g_i^{-1},\,1\leq i\leq N.$\\
One gets the finite dimensional Hopf algebras in \cite{AS} from this construction if one considers only finite abelian groups, linkings that link different connection components of the given finite Dynkin diagram, and if one introduces the extra root vector relations \cite[(1.8)]{AS}.

\section{The structure of linkable Dynkin diagrams}
In this section we want to explain the structure of linkable Dynkin diagrams for which there exist linkable braiding matrices.
This will allow us to construct a large number of examples for these new Hopf algebras. To get a nice result however, we slightly specialize some of  our earlier definitions. We will discuss generalisations in the next chapter.
\begin{quote}From now on all linkable Dynkin diagrams are assumed to be link-connected, i.e. when viewed as a graph they are connected. Furthermore we will restrict our considerations to diagrams where two vertices are linkable only if they lie in different connection components of the original diagram.
\end{quote}
For two vertices $i,j$ of the Dynkin diagram with $a_{ij}\neq0$ the symmetry of (\ref{cartan}) implies \begin{equation}\label{normdiag}b_{ii}^{a_{ij}}=b_{jj}^{a_{ji}}.\end{equation}For $i\lkb j$ we have $a_{ij}=0$, as we required the vertices to lie in different connection components. Using (\ref{link}) and (\ref{cartan}) alternately, we arrive at \begin{equation}\label{diaglink}b_{ii}=b_{ij}^{-1}=b_{ji}=b_{jj}^{-1}.\end{equation}
\subsection{The finite case}
First we will only consider Dynkin diagrams of finite type, i.e. the corresponding Lie algebras are finite dimensional.
In order to get interesting applications in regard of \cite{AS} we further require that a linkable braiding matrix has the following property:\\

\refstepcounter{equation}
\hfill\begin{parbox}{.8\textwidth}{\emph{The order of the diagonal elements $b_{ij}$ is greater than 2 and not divisible by 3 if the linkable Dynkin diagram contains a component of type $G_2$.}}
\end{parbox}\hfill  (\theequation)\\ \label{order}

The first properties are presented in a lemma, which is essentially Lemma 5.6. in \cite{AS}. However, we formulate it on the level of the braiding matrix.

\beer{lemma}{\label{aslemma}We are given a linkable Dynkin diagram $D$ and a corresponding linkable braiding matrix ${\bf b}.$ Suppose that the vertices $i$ and $j$ are linkable to $k$ and $l$, respectively. Then $a_{ij}=a_{kl}.$}
\prf{If $a_{il}\neq 0$ or $a_{jk}\neq{0}$ then we immediately get $a_{ij}=a_{kl}=0$, because linkable vertices must lie in different connection components of $D.$ So we now take $a_{il}=a_{jk}=0.$
Without loss of generality we assume $a_{ij}\le a_{kl}.$
Using (\ref{cartan}) and (\ref{link}) alternately, we get
$$b_{ii}^{a_{ij}}=b_{ij}b_{ji}=b_{il}^{-1}b_{jk}^{-1}=b_{li}b_{kj}=b_{lk}^{-1}b_{kl}^{-1}=b_{kk}^{-a_{kl}}=b_{ii}^{a_{kl}}.$$
In the last step we used (\ref{diaglink}). Hence $a_{ij}=a_{kl}$ modulo the order of $b_{ii}.$  As $b_{ii}\neq\pm 1$ we either get $a_{ij}=a_{kl}$ or that the order of $b_{ii}$ is 3 and $a_{ij}=-3, a_{kl}=0.$ But in the last case $i$ and $j$ form a $G_2$ component.
So $b_{ii}=3$ is a contradiction to the assumption on the order of the diagonal elements.
}
Before we can state our result on the structure of linkable Dynkin diagrams that admit a corresponding braiding matrix with the above properties, we have to introduce some terminology.
\defi{For every cycle\footnote{A cycle is a closed, non self-intersecting path in the diagram.} $c$ in $D$ we choose an orientation and denote by the \emph{weight} $w_c$ the absolute value of the difference of the numbers of double edges in that cycle with the arrow pointing with the orientation and against it. The \emph{length} $l_c$ of the cycle is defined to be the number of dotted edges in that cycle.\\
The \emph{genus} $g_c$ of the cycle is now defined by the following formula:
\begin{equation}
g_c:=2^{w_c}-(-1)^{l_c}.
\end{equation}}
In preparation for some technicalities in the second part of the proof of our result we need the following concept.
\defi{
For two vertices $i$ and $j$ of $D$ we define for every directed path $P$ from $j$ to $i$ a number $h^i_j(P)\geq 0$, called the \emph{height of $i$ over $j$ along $P$}, by the following algorithm.\\
First we set $h=0$. Then we follow the path $P$  starting at $j$. At every vertex we get to, we 
\begin{quote}\begin{center}{\footnotesize decrease the value of $h$ by 1\\ increase it by 1\\ or leave it unchanged,}\end{center}\end{quote}
 depending on if the edge we just passed was a double edge pointing  
\begin{quote}\begin{center}{\footnotesize with the orientation of $P$\\ against it\\ or was not a double edge.}\end{center}\end{quote} 
The only exception is that the value of $h$ is not decreased when it is 0. $h^i_j(P)$ is then set to be the value of $h$ after we followed through the whole Path $P$ arriving at $i.$\\
For a cycle $c$ we define the \emph{natural orientation} to be the one where the number of double edges in $c$ pointing with this orientation is not less than the number of double edges pointing against it\footnote{If the weight $w_c=0$ then the natural orientation is ambiguous. In that case we choose one of the possible two orientations. This will not lead to any problems.}.\\
For every vertex $i$ of $c$ we define the \emph{absolute height} $h_i(c)\geq 0$ to be the height of $i$ over itself along $c$ following its natural orientation. A vertex of absolute height $0$ in a cycle of genus $g_c>0$ is called a \emph{Level 0 vertex.}}

This seems to be the right point to illustrate all the notions in an example. We consider the following Dynkin diagram where the vertices are supposed to be linkable in the indicated way:
\abild{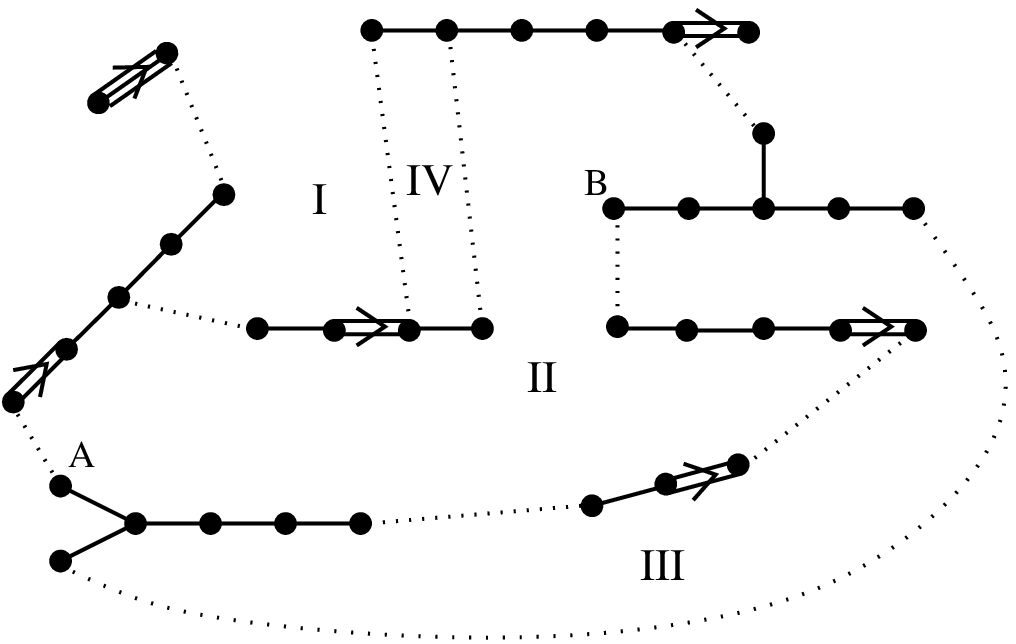}
For the four cycles denoted by I, II, III and IV (where I is the outside cycle) the values for $w$, $l$ and $g$ are given in this table:
\begin{center}\begin{tabular}[b]{c|c c c}
&$w_c$&$l_c$&$g_c$\\
\hline
I&2&5&5\\
II&2&7&5\\
III&0&4&0\\
IV&0&2&0\\
\end{tabular}.\end{center}

The natural orientation of cycles I and II is clockwise, whereas the natural orientation in cycles III and IV is ambiguous. The vertex indicated by the letter ``A'' is a vertex of absolute height 1 in cycle II, but a Level 0 vertex for cycle I. And Vertex ``B'' is a Level 0 vertex for cycle II but a vertex of absolute height 1 for cycle III independently of the natural orientation chosen for that cycle.

We are now able to come to our main result.
\beer{thm}{
We are given a link-connected linkable Dynkin diagram $D$ and explicitly exclude the case $G_2\times G_2$. It will be treated later.\\
There exists a linkable braiding matrix of $D$-Cartan type, iff
\begin{enumerate}
\item\label{cd1} In components of type $G_2$ not both vertices are linkable to other vertices.
\item\label{cd2} $D$ does not contain any induced subgraphs\footnote{An induced subgraph consists of a subset of the original vertices and all the corresponding edges.} of the form:
\abild{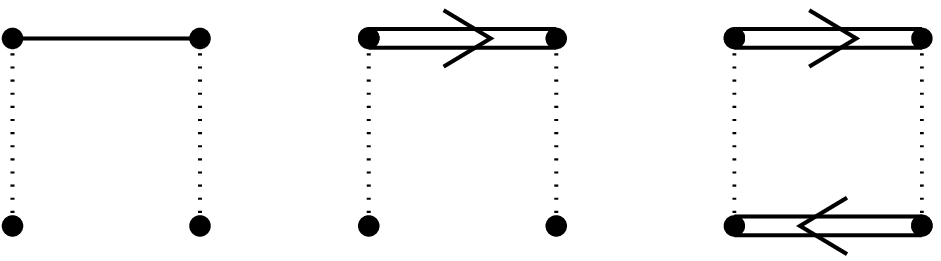}
\item\label{cd3} One of the following conditions is true: 
\begin{itemize}
\item[] $D$ contains no cycles or the genera of all cycles are zero.
\item[or] $D$ does not contain a component of type $G_2$ and there is a common divisor $d>2$ of all cycle genera and the field $k$ contains a primitive $\pth d$ root of unity.
\item[or] $D$ does contain a component of type $G_2$ and there is a common divisor $d>2$ of all cycle genera, $d$ is not divisible by 3 and the field $k$ contains a primitive $\pth d$ root of unity.
\end{itemize}
\end{enumerate}}
\prf{
We first prove the ``if'' part, i.e. we assume conditions \ref{cd1}.-\ref{cd3}.\\

We will construct the braiding matrix explicitly and show that it  fulfills the required identities.\\
The main observation is that once an element of the diagonal is chosen, the other diagonal elements are determined (up to possible signs) by (\ref{normdiag}) and (\ref{diaglink}).\\
We take $d>2$ as given by condition \ref{cd3}. In the case that $D$ contains no cycles or the genera of all cycles are zero we set $d>2$ to be a prime, such that $k$ contains a primitive $\pth d$ root of unity. This is possible by the general assumption on the field $k$.  We note in particular that $d$ is always odd.\\
Now we choose a vertex $i$ 
and set $b_{ii}:=q$, where $q$ is a primitive $\pth d$ root of unity. As the Dynkin diagram $D$ is link-connected, we can choose for every vertex $j\neq i$ a path\footnote{Again we demand that a path does not include a vertex more than once.} $P_{ij}$ connecting $i$ and $j$, which we denote by the sequence of its vertices $(i=p_0,p_1,\dots,p_t=j).$ For every such path $P_{ij}$ we now define the $b_{p_kp_k},\; k=1,\dots,t,$ recursively: 
\begin{equation}\label{recurs}b_{p_{k+1}p_{k+1}}=\begin{cases} b_{p_kp_k}^{-1},&\text{if }a_{p_kp_{k+1}}=0\\
                                   b_{p_kp_k}^{a_{p_kp_{k+1}}/a_{p_{k+1}p_k}},&\text{otherwise.}\end{cases}\end{equation}
When taking the square root we choose it to be again a $\pth d$ root of unity. As $d$ is odd, this picks exactly one of the two possible roots. When taking third roots we pick them to be the unique $\pth d$ root of unity as well. This is always possible, as we only have to take third roots when we are dealing with vertices of components of type $G_2.$ Then however, $d$ is not divisible by 3.

We now show that this process is well defined.\\
The only problems could arise when a vertex can be reached from the starting vertex by two different paths, i.e. when $D$ contains a cycle. Suppose we are given two different paths $(i=p_0,p_1,\dots,p_t=j)$ and $(i=q_0,q_1,\dots,q_u=j).$ Let $n$ be the smallest integer with $p_{n+1}\neq q_{n+1}$ and $m_1>n$ the smallest integer, such that there is a $m_2>n$ with $p_{m_1}=q_{m_2}$. Then $c=(p_n,p_{n+1},\dots,p_{m_1}=q_{m_2},q_{m_2-1},\dots,q_n=p_n)$ is a cycle. It is now sufficient to show that the recursive procedure (\ref{recurs}) for the paths $P_1:=P_{p_np_{m_1}}$ and $P_2:=P_{q_nq_{m_2}}$ leads to the same value 
\begin{equation}\label{equal}b_{p_{m_1}p_{m_1}}=b_{q_{m_2}q_{m_2}}.\end{equation}
As triple edges are not part of cycles, we easily obtain a closed formula for the desired values:
\begin{equation}b_{p_{m_1}p_{m_1}}=Q^{(-1)^{l_1}2^{w_1}}\qquad \text{and}\qquad b_{q_{m_2}q_{m_2}}=Q^{(-1)^{l_2}2^{w_2}}.\end{equation}
Here $Q:=b_{p_np_n}$, $l_i\ge 0$ denotes the number of dotted edges in the path $P_i$ and $w_i\in\Z$ is the difference of the numbers of double edges in $P_i$ that have the arrow pointing with the path's orientation and against it. Without loss of generality we assume $w_1\ge w_2$ and have $w_c=w_1-w_2$ and $l_{c}=l_1+l_2$.\\
As $Q$ is a $\pth d$ root of unity and $d$ divides all cycle genera we get $Q^{g_c}=1$ or $$Q^{2^{(w_1-w_2)}}=Q^{(-1)^{(l_1+l_2)}}.$$
Taking both sides to the power of $(-1)^{l_1}2^{w_2},$ we arrive at (\ref{equal}). Here we would like to remind the reader that all values are $\pth d$ roots of unity and hence there is no ambiguity regarding signs.\\
The so set diagonal entries of the braiding matrix fulfill the requirements for their orders, because $d$ is odd and not divisible by 3 when there are components of type $G_2$ in $D.$\\

We now give the remaining entries for the linkable braiding matrix, i.e. we set $b_{ij}$ for $i\neq j.$ For this we divide the set $\{(i,j):i\neq j\}$ of pairs of vertices into 4 classes:
\begin{description}
\item[None of the two vertices
is linkable to some other one.]
We set $$b_{ji}:=z  \qquad b_{ij}:=b_{ii}^{a_{ij}}z^{-1}.$$
\item[The two vertices
are linkable to each other.]
We set \begin{align*}b_{ij}&:=b_{ii}^{-1}\\
                        b_{ji}&:=b_{jj}^{-1}.\end{align*}   
\item[Only one of the two vertices
is linkable to some other vertex.] We assume $i$ is linkable to $k$. We set
\begin{align*} b_{ji}&:=z  &\qquad b_{ij}&:=b_{ii}^{a_{ij}}z^{-1}\\
                  b_{jk}&:=z^{-1} &\qquad b_{kj}&:=b_{kk}^{a_{kj}}z \end{align*}
\item[Both vertices
are linkable to some other vertices.]
We assume $i$ is linkable to $k$ and $j$ is linkable to $l.$ For $i$ and $k$ to be linkable we can not have $a_{ij}\neq 0$ and $a_{jk}\neq 0.$ So after a possible renaming of the indices $i$ and $k$ we can assume that $a_{jk}=0.$ By the same reasoning we take $a_{il}=0.$ Now we set
\begin{align*} b_{ji}&:=b_{kj}:=z &\qquad b_{ij}&:=b_{li}:=b_{ii}^{a_{ij}}z^{-1}\\
                  b_{jk}&:=b_{kl}:=z^{-1} &\qquad b_{il}&:=b_{lk}:=b_{ii}^{-a_{ij}}z.\end{align*}
\end{description}
In all the cases $z\neq 0$ can be chosen freely from the field $k$ and can be different for every class and pair of vertices.\\
We would like to point out that all pairs of indices fall into one of those classes and that there are no overlapping cases, i.e. each off-diagonal element is only set in one of these.\\
This way we have explicitly constructed the matrix ${\bf b}=(b_{ij})$. We are left to show that (\ref{cartan}) and (\ref{link}) are fulfilled. For the diagonal entries this has been done already. For the entries being set in the first three classes it is immediately clear from the definition.\\
In the last class only the relation $b_{kl}b_{lk}=b_{kk}^{a_{kl}}$ must still be checked.\\
We note that neither $i$ and $j$ nor $k$ and $l$ can form a component of type $G_2$, as this would contradict part \ref{cd1}. of the assumption.\\ From the construction we get $b_{kl}b_{lk}=b_{ii}^{-a_{ij}}.$ As vertex $i$ is linkable to $k$ we know $b_{kk}=b_{ii}^{-1}.$ We will show that $a_{kl}=a_{ij}.$\\ 
If $a_{ij}=0$ then we immediately get $a_{kl}=0$, because $a_{kl}a_{lk}=1$ or 2 is not permitted by part \ref{cd2}. of the assumption. Analogously we get the result if we assume $a_{kl}=0.$\\
The case where the 4 indices form a sub-diagram of the kind
\abild{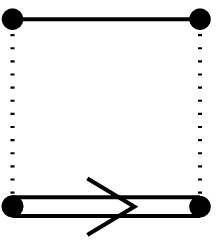}
is excluded, as the cycle genus for this diagram is $1=2^1-(-1)^2.$ So the only other possible diagrams these four vertices can form are
\abild{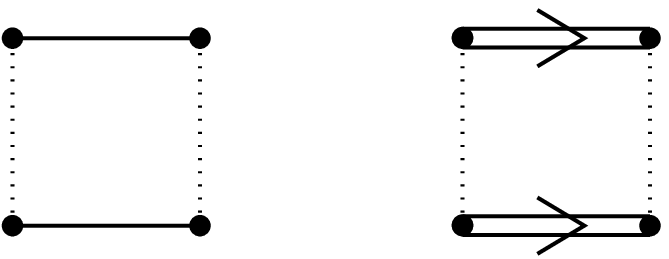}
which all have $a_{ij}=a_{kl}.$ We would like to note that the cycle genus for the last diagram is $0=2^0-(-1)^2.$\\
This concludes the ``if'' part of the proof.\\

Before we come to the ``only if'' part we prove a lemma to enable us to deal with some arising technicalities.\\
\beer{lemma}{\label{tech}
We are given a linkable Dynkin diagram $D$.
\begin{enumerate}
\item In every cycle $c$ of $D$ with $g_c>0$ there exists a Level 0 vertex.
\item \label{simplecyc}Given a linkable braiding matrix ${\bf b}$ of $D$-Cartan type we have for every Level 0 vertex $i$ of $c$: \quad $b_{ii}^{g_c}=1.$
\item \label{allcyc}Let $G$ be the greatest common divisor of all cycle genera. If there are no cycles with cycle genus 1 or 2 we have for every Level 0 vertex $i$:\quad $b_{ii}^G=1.$
\end{enumerate}}
\prf[of Lemma]{
\begin{enumerate}
\item Pick a vertex $i$ in $c$
and calculate $h_i(c)$. If $h_i(c)>0$ then take the vertex $j$, where (in the recursive definition) $h$ was 0 for the last time.\\ $j$ is then a Level 0 vertex, because:
\begin{itemize}
\item In the recursive calculation of $h_j(c)$ the value of $h$ is positive at least until we pass vertex $i$. (Choice of $j$ and $i$ is not a Level 0 vertex.)
\item Assume $h$ stays positive until it reaches $j$ again, i.e. $j$ is not Level 0. Then however, the number of double edges along the cycle $c$ pointing against the natural orientation is greater than the number of double edges pointing with it, which is a contradiction to the definition of natural orientation.
\item So there is a vertex $k$ between $i$ and $j$ where $h$ becomes 0. But this means that the value of $h$ at $k$ in the calculation of $h_i(c)$ was 0 as well. (This value can not be bigger, having started out smaller at vertex $i.$)
\item Now, the calculation of $h$ from $k$ until $j$ is the same as for $h_i(c).$ So $h=0$ when it reaches $j.$ 
\end{itemize}
   
\item Set $Q:=b_{ii}.$ Following the cycle in its natural orientation starting at $i$ and using (\ref{normdiag}) and (\ref{diaglink}) we arrive at
$$Q=Q^{(-1)^{l_c}2^{w_c}}.$$
That there are no extra signs from possible square roots in the above formula is ensured by the assumption that $i$ is of Level 0. Raising both sides to the $\pth{{(-1)^{l_c}}}$ power and dividing by the new left hand side we get $1=Q^{g_c}.$
\item If there are no cycles $c$ with $g_c>0$ there is nothing to show. When there is only one such cycle $c$ we have $G=g_c$ and part \ref{simplecyc}. of this Lemma establishes the claim.\\
We take now two cycles $c_1$ and $c_2$ with $g_{c_1}\ge g_{c_2}>2$ and set $g:=\gcd(g_{c_1},g_{c_2}).$ 
We pick in $c_1$ a vertex $i$ of Level 0 and a vertex $j$ of Level 0 in $c_2.$ Following a path from $i$ to $j$ and applying (\ref{normdiag}) and (\ref{diaglink}) appropriately we get $b_{jj}^{2^{w_1}}=b_{ii}^{(-1)^l2^{w_2}}$ for some values $w_1,w_2$ and $l.$ $b_{ii}$ is a $\pth{g_{c_1}}$ root of unity (see previous part of this Lemma) and so $b_{jj}$ is also a $\pth{g_{c_1}}$ root of unity. However, using the previous part again, $b_{jj}$ must be a $\pth{g_{c_2}}$ root of unity. As the cycle genera are not divisible by 2 we conclude that $b_{jj}$ is a $\pth g$ root of unity.\\
Repeating the argument for all the other cycles with cycle genus bigger than 2 we conclude that for every Level 0 vertex $i$ the corresponding diagonal entry $b_{ii}$ of the braiding matrix is a $\pth G$ root of unity. 
\end{enumerate}}

Now we finish the proof of the theorem. We are given a linkable braiding matrix $\mathbb b$ and a corresponding linkable Dynkin diagram $D$ and set $G$ to be the greatest common divisor of all cycle genera. $G:=0$ if all cycle genera are 0.\\
Assume now:
\begin{itemize}
\item there is a $G_2$ component with both vertices $i$ and $j$ linkable to other vertices $k$ and $l$, respectively.\\
From Lemma \ref{aslemma} we immediately get $a_{kl}=a_{ij}$ and $a_{lk}=a_{ji}.$ So $k$ and $l$ form another $G_2$ component. Thus the given diagram must be $G_2\times G_2,$ which we don't want to consider here.
\item one induced subgraph is of the kind as in condition \ref{cd2}.\\
This contradicts Lemma \ref{aslemma}.
\item there is a cycle $c$ with cycle genus $g_c=1$ or 2.\\
Using Lemma \ref{tech} \ref{simplecyc}. we get that a diagonal entry of the braiding matrix must be 1 or $-1$.
\item $G=1$ and there is no cycle with genus 1 or 2.\\
Lemma \ref{tech} \ref{allcyc}. shows that then $b_{ii}=1$ for all Level 0 vertices $i$, a contradiction to our assumption (\ref{order}) about the order.
\item $G>2$ and the base field $k$ does not contain a primitive $\pth d$ root of unity for any $d>2$ dividing $G$.\\
Lemma \ref{tech} \ref{allcyc}. shows that the field $k$ must contain a $\pth G$ root of unity. This means there must be a $d$ dividing $G$, such that $k$ contains a primitive $\pth d$ root of unity. By the order assumption (\ref{order}) 
there is even a $d$ bigger than 2.
\item $G>2$, there is a component of type $G_2$ and the only $d>2$ that divide $G$, such that $k$ contains a primitive $\pth d$ root of unity, are divisible by 3.\\
There is a cycle $c$ with $g_c>2.$ According to Lemma \ref{tech} there is a Level 0 vertex $i$ in $c$ with $b_{ii}$ a primitive $\pth d$ root of unity, where $d>2$ and $d$ divides $G$. From the assumption we get that $d$ must be divisible by 3, so the order of $b_{ii}$ is divisible by 3. This contradicts (\ref{order}).
 \end{itemize}}

\subsection{The affine case}
Now we turn to the affine case, i.e. we consider Dynkin diagrams that are unions of diagrams of finite and affine type. In order to get a similar result to the previous one we have to consider an even more specialized style of braiding matrix.\\

\refstepcounter{equation}
\hfill\begin{parbox}{.8\textwidth}{\emph{We require that the order of all diagonal entries is the same and equal to a prime bigger than 3. This kind of braiding matrix we will call \emph{homogeneous}.}}
\end{parbox}\hfill  (\theequation)\\ \label{afforder}

We have the analog of Lemma \ref{aslemma} for this situation.
\beer{lemma}{\label{asaff}We are given a linkable Dynkin diagram $D$ and a corresponding homogeneous linkable braiding matrix ${\bf b}.$ Suppose that the vertices $i$ and $j$ are linkable to $k$ and $l$, respectively. Then $a_{ij}=a_{kl}.$}

The proof is the same as before and the last conclusion is straight forward, as $b_{ii}$ has to have at least order 5, according to the assumptions.

As before we define the notions of weight and natural orientation for every cycle $c$. This time, however, we do this as well for triple edges. The former notions will now be denoted by \emph{natural 2-orientation} and $w^2_c.$ The corresponding ones for the triple edges by \emph{natural 3-orientation} and $w^3_c.$ The length of the cycle is exactly as before.\\
If the two natural orientations coincide we define the genus of this cycle as $$g_c:=3^{w^3_c}2^{w^2_c}-(-1)^{l_c}.$$
In the other case we take $$g_c:=|3^{w^3_c}-2^{w^2_c}(-1)^{l_c}|.$$

Now the theorem can be formulated in the same spirit.

\beer{thm}{
We are given a link-connected linkable affine Dynkin diagram $D$ and explicitly exclude the cases $A_1^{(1)}\times A_1^{(1)}$ and $A_2^{(2)}\times A_2^{(2)}$. These will be treated later.\\
There exists a homogeneous linkable braiding matrix of $D$-Cartan type, iff

\begin{enumerate}
\item\label{dc1} In components of type $A_1^{(1)}$ and $A_2^{(2)}$ not both vertices are linkable to other vertices.
\item\label{dc2} $D$ does not contain any induced subgraphs of the form:
\abild{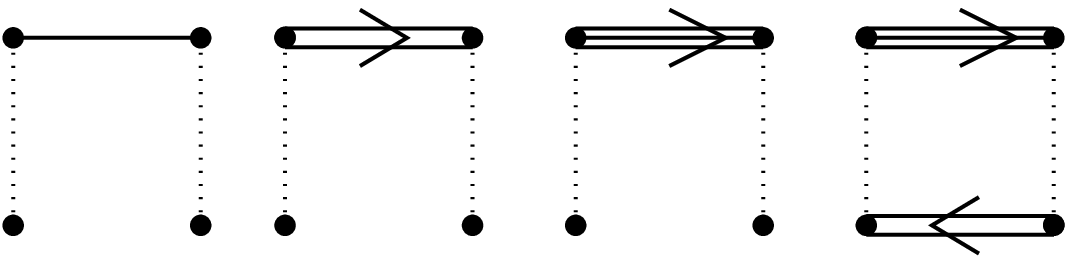}
\item\label{dc3} There is a prime $p>3$ that does divide all cycle genera and the field $k$ contains a $\pth p$ root of unity.
\end{enumerate}}

\prf{(sketch)
We start again with the ``if'' part. So suppose all conditions are fulfilled.\\
We pick any vertex $i$ and set $b_{ii}$ to be a $\pth p$ root of unity. As before, we can now set all other diagonal elements of the braiding matrix inductively. Dotted, single, double and triple edges are treated as in the finite case, the edge of type $A_1^{(1)}$ is treated as a single edge and when we pass quadruple edges we take the $\pth 4$ power or root according to if we go along with the arrow or in the opposite direction. This time we only require that all diagonal values are $\pth p$ roots of unity and as $p>3$ is a prime this determines the values uniquely. The independence from the paths chosen in this procedure is ensured again by the condition on the genera and can be shown by elaborating on the technique used in the finite case.\\
The off-diagonal elements are set completely in the same way as in the finite case, so we don't repeat the arguments here.\\

The ``only if'' part follows the same strategy as before as well. We assume to be given a linkable Dynkin diagram $D$ and a homogeneous linkable braiding matrix of $D$-Cartan type.
Suppose condition \ref{dc1}. is not fulfilled. Then Lemma \ref{asaff} immediately gives that the diagram is $A_1^{(1)}\times A_1^{(1)}$ or $A_2^{(2)}\times A_2^{(2)}.$\\
The negation of condition \ref{dc2}. contradicts Lemma \ref{asaff}.\\
Suppose now that the biggest prime $p$ dividing all cycle genera is smaller than 5 or that the greatest common divisor $G$ of all cycle genera is 1. 
Here the analog of Lemma \ref{tech} \ref{allcyc}. gives the contradiction. One does not need to use the concept of Level 0 vertices, however, as all diagonal entries have to have the same order. Still the result that every $b_{ii}$ is a $\pth G$ root of unity can be deduced by the same reasoning.} 
\subsection{The excluded cases}
Suppose the given diagram is one of the above excluded ones, i.e. $G_2\times G_2,$ $A_1^{(1)}\times A_1^{(1)}$ or $A_2^{(2)}\times A_2^{(2)}.$\\
If the linkable diagram is not excluded by Lemma \ref{asaff} then the following linkable braiding matrices do exist for any of the possible linkings
\begin{equation*} \begin{pmatrix} q&z&q&z^{-1}q^{-m}\\z^{-1}q^{-m}&q^n&z&q^n\\q^{-1}&z^{-1}&q^{-1}&zq^m\\zq^m&q^{-n}&z^{-1}&q^{-n}\end{pmatrix}. \end{equation*}
Here $n=m=3$; $n=1,m=2$ or $n=m=4$ for the first, second and last diagram, respectively.

\subsection{Examples}
If we start with two copies of a Dynkin diagram and link corresponding vertices, all cycles have genus 0, so there are no obstructions. This shows once more that braiding matrices necessary for the construction of the quantum universal enveloping algebras $\uqg$ exist.\\
If we take $n$ copies of $A_3$ and link them into a circle \cite[Example 5.13.]{AS}, then we have one cycle of weight 0 and length $n$. So the genus is 0 for $n$ even and there are no restrictions. For $n$ odd, however, the genus is 2, so under the conditions we imposed there are no corresponding linkable braiding matrices.\\
For the case of $n$ copies of $B_3$ linked into a circle the genus is $2^n-(-1)^n.$\\ 
In the example which we used to demonstrate the various notions, we found that all cycle genera were divisible by 5. So we know that all linkable braiding matrices of that type have diagonal elements, whose order is divisible by 5. This gives a limitation on the groups over which these matrices can be realized.

\section{Generalisations}
Here we want to discuss some possible generalisations to restrictions imposed in the last chapter. We only considered affine and finite Cartan matrices, as we do not know of any further special classes within generalized Cartan matrices. Apart from that, the combinatorics involved in the classification become more involved with increasing values of the entries of the Cartan matrix. In order to get a nice presentation, one would have to impose even more restrictions on the braiding matrix.
\subsection{The order of the diagonal elements}
We would like to make some comments on why we made the various restrictions above on the orders of the diagonal elements.\\
First, if $b_{ii}=1$ then $g_i$ commutes with $a_i,$ a not very interesting fact. Moreover \cite[Lemma 3.1.]{ASp3} shows that this case can not emerge when dealing with finite dimensional Hopf algebras.\\
When the order of $b_{ii}$ can be two, we don't get Lemma \ref{aslemma}. So we would have to deal with a much more difficult structure of possible diagrams.\\ Actually, if one neglects the sub-diagram
\abild{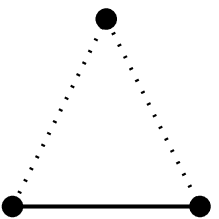}
then Lemma \ref{aslemma} gives us that a vertex can be linkable to at most one other vertex. One just takes $i=j$ and finds $a_{kl}=2$ or $a_{kl}=a_{lk}=-1$ and the order of $b_{ii}$ is 3. If the order of $b_{ii}$ is 2, then one could consider a great number of diagrams with vertex $i$ being linkable to more than one vertex, thus making a classification even harder.\\
To simplify the presentation of the theorems we excluded orders divisible by 3 when there are components of type $G_2.$ Without this limitation, a much more thorough examination of the diagrams (with heavy use of the above defined heights) is needed to establish a necessary condition for the existence of a braiding matrix. The problems come from trying to avoid diagonal elements of order 1 and 2.\\
For the affine diagrams this problem is even more severe, and the easiest way to avoid it is to consider only homogeneous braiding matrices with prime order of diagonal elements greater than 3. This way, there are no extra difficulties stemming from triple arrows either.

\subsection{Linkings within one component}
We first consider the case where two vertices $i$ and $j$ within the same connection component are linkable, but not neighbouring, i.e. $\a =0.$ In that case we can actually still apply last sections's considerations. In Dynkin diagrams of affine and finite type the possible values of the genera of cycles formed by this special linking can be calculated easily. We find as possible values
\begin{equation}
g_c=\left\{\begin{array}{l}2\\3\\4\\5\end{array}\right.\text{if the diagram between $i$ and $j$ }\begin{array}{l}\text{is of type $A,$ $C^{(1)}$ or $D^{(2)}$}\\\text{has one double edge}\\\text{has one triple edge}\\\text{is of type $A_{2k}^{(2)},\;k\ge 2.$}\end{array}
\end{equation}
If $\a\neq0$ and we link $i$ and $j$ we get from (\ref{link}) by setting $k$ first to $i$ then $j$ and exchanging $i$ and $j$ the following identity
\begin{equation} b_{ii}^{\a a_{ji}-\a-a_{ji}}=1.\end{equation}
For the possible sub-diagrams
\abild{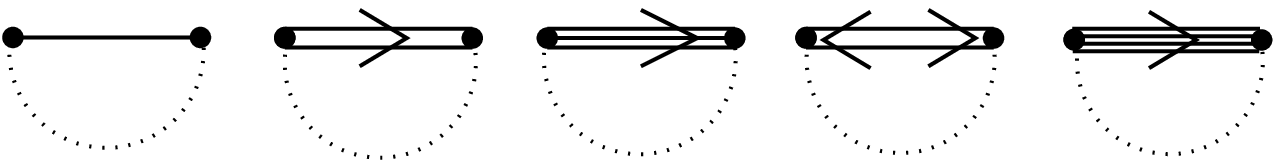}
the condition gives us that the order of $b_{ii}$ must divide 3, 5, 7, 8 or 9, respectively.\\
So we see that admitting self-linkings restricts the possible braiding matrices drastically, especially if we are interested in homogeneous ones.\\
The values of 3 and 5 we get, when we link the two vertices in diagrams $A_2$ and $B_2$, respectively, are as well the problematic ones when one tries to find all liftings of the Nichols algebras of this type, cf. \cite{ASa2,BDR}. 
\subsection{Link-disconnected diagrams}
The results on link-connected diagrams can easily be extended to arbitrary diagrams. For each link-connected component the considerations can be carried out and a possible braiding matrix constructed. The direct sum of these matrices is then a braiding matrix for the whole diagram. If we are interested only in homogeneous braiding matrices, then we have to check that we can choose the orders of the diagonal elements in all matrices (corresponding to the various link-connected components) to be the same.

\section{Group realization} 

We want to close this article with a discussion of the role of the group $\G$. First of all it should be noted that any linkable $(s\times s)$ braiding matrix $(b_{ij})$ can be realized over $\Z^s$ simply by taking the elements of the canonical basis of this group as the $g_i$ and defining the characters by $\chi_j(g_i):=b_{ij}.$ This can actually be extended to groups $(\Z/p)^s$ as long as $p$ is divisible by the orders of the elements of the braiding matrix.\\ 
If, however, we are given a group $\G$ it is not clear to us as to which braiding matrices can be realized over this group. In \cite[chap. 8]{AS-cartan} some answers in this direction were found.\\
We would like to present here some details for the group $\G=(\Z/p)^2,\,p>3.$ According to the considerations in \cite[Prop. 8.3.]{AS-cartan} the biggest finite Dynkin diagram (for which a braiding matrix exists that is realizable over $\G$) has four vertices unless $p=5$ and the diagram is $A_4\times A_1.$ A computer check\footnote{A program explicitly constructs the braiding matrices realizable over the group.} seems to indicate that this bound is the only limitation if one neglects the diagram $A_4.$ This raises the hope that for $(\Z/p)^s$ with $s>2$ and $p>3$ all finite Dynkin diagrams within the bounds derived in the above mentioned paper admit a braiding matrix realizable over this group. We hope to contribute to this question in a later work.\\
For the Dynkin diagram $A_4$ we would like to present the calculation that reveals the limitations in that case.\\
We take a canonical basis $a$ and $b$ of $\G$ as $g_1$ and $g_2$ and the elements $a^nb^m$ and $a^kb^l$ for some $0\leq n,m,k,l< p$ as $g_3$ and $g_4$. Suppose now, that we have a braiding matrix ${\bf b}$ of type $A_4$ that can be realized with these $g_i$ and some characters $\chi_i$. This means in particular that $p$ must be divisible by the orders of all entries $b_{ij}$.
As (\ref{cartan}) must be fulfilled we get the following relations modulo $p:$
\begin{align}
\label{magic}n^2-nm+m^2+m+1&=0\\
\label{pure}k^2-kl+l^2+1&=0\\
\label{mix}k(m-2n)+l(n-2m-1)+1&=0
\end{align}
The first equation (\ref{magic}) has for every prime $p=6z\pm 1$ exactly $6z$ solutions $(n,m)$.
We assume the above relations to hold and now investigate the following cases:
\begin{description}
\item[Case 1: $m=2n$.] We immediately get from (\ref{magic}) $n=\frac{-1\pm\sqrt{-2}}3$ and from (\ref{mix}) $l=\pm\frac 1{\sqrt{-2}}.$ Using this in (\ref{pure}) we get $k=\frac{\pm 1\pm\sqrt 5}{2\sqrt{-2}}.$
\item[Case 2: $n=2m+1$.] With the same reasoning we get this time $m=\frac{-2\pm\sqrt{-2}}3$, $k=\pm\frac 1{\sqrt{-2}}$ and because of the symmetry of (\ref{pure}) $l=\frac{\pm 1\pm\sqrt 5}{2\sqrt{-2}}.$
\item[Case 3: $m\neq 2n,\;n\neq 2m+1$.] Solving (\ref{mix}) for $k$ and substituting it into (\ref{pure}) we find by using (\ref{magic}) that $l=-\frac 12-\frac{3m}4\pm\sqrt 5\frac{m-2n}4$ and $k=\frac{3m(n-2m-1)-2}{4(m-2n)}\mp\sqrt 5\frac{n-2m-1}4.$
\end{description}
Thus we see that for the system to have a solution, 5 must be a square modulo $p$. According to the reciprocity law this happens for $$p\equiv 1 \text{ or 3 }\; \text{ mod 10}\qquad\text{or}\quad p=5.$$ This condition is as well sufficient as (\ref{magic}) always has a solution.

\end{document}